\documentclass[3p,,preprint,12pt]{elsarticle}
\makeatletter\if@twocolumn\PassOptionsToPackage{switch}{lineno}\else\fi\makeatother

\usepackage{tabulary,xcolor}
\usepackage{amsfonts,amsmath,amssymb}
\usepackage[T1]{fontenc}
\makeatletter
\let\save@ps@pprintTitle\ps@pprintTitle
\def\ps@pprintTitle{\save@ps@pprintTitle\gdef\@oddfoot{\footnotesize\itshape \null\hfill\today}}
\def\hlinewd#1{%
  \noalign{\ifnum0=`}\fi\hrule \@height #1%
  \futurelet\reserved@a\@xhline}

\AtBeginDocument{\ifNAT@numbers \biboptions{sort&compress}\fi}
\makeatother

\usepackage{ifluatex}
\ifluatex
\usepackage{fontspec}
\defaultfontfeatures{Ligatures=TeX}
\usepackage[]{unicode-math}
\unimathsetup{math-style=TeX}
\else 
\usepackage[utf8]{inputenc}
\fi 
\ifluatex\else\usepackage{stmaryrd}\fi

\usepackage{url,multirow,morefloats,floatflt,cancel,tfrupee}
\makeatletter

\AtBeginDocument{\@ifpackageloaded{textcomp}{}{\usepackage{textcomp}}}
\makeatother
\usepackage{colortbl}
\usepackage{xcolor}
\usepackage{pifont}
\usepackage[nointegrals]{wasysym}
\makeatletter

\def\mcWidth#1{\csname TY@F#1\endcsname+\tabcolsep}

\def\cAlignHack{\rightskip\@flushglue\leftskip\@flushglue\parindent\z@\parfillskip\z@skip}
\def\rAlignHack{\rightskip\z@skip\leftskip\@flushglue \parindent\z@\parfillskip\z@skip}

\@ifundefined{etal}{}{}

\usepackage{ifxetex}
\ifxetex\else\if@twocolumn\@ifpackageloaded{stfloats}{}{\usepackage{dblfloatfix}}\fi\fi

\AtBeginDocument{
\expandafter\ifx\csname eqalign\endcsname\relax
\def\eqalign#1{\null\vcenter{\def\\{\cr}\openup\jot\m@th
  \ialign{\strut$\displaystyle{##}$\hfil&$\displaystyle{{}##}$\hfil
      \crcr#1\crcr}}\,}
\fi
}

\AtBeginDocument{%
  \@ifpackageloaded{endfloat}%
   {\renewcommand\efloat@iwrite[1]{\immediate\expandafter\protected@write\csname efloat@post#1\endcsname{}}}{\newif\ifefloat@tables}%
}%

\def\BreakURLText#1{\@tfor\brk@tempa:=#1\do{\brk@tempa\hskip0pt}}
\let\lt=<
\let\gt=>
\def\processVert{\ifmmode|\else\textbar\fi}

\@ifundefined{subparagraph}{
\def\subparagraph{\@startsection{paragraph}{5}{2\parindent}{0ex plus 0.1ex minus 0.1ex}%
{0ex}{\normalfont\small\itshape}}%
}{}

\newcommand\role[1]{\unskip}
\newcommand\aucollab[1]{\unskip}
  
\@ifundefined{tsGraphicsScaleX}{\gdef\tsGraphicsScaleX{1}}{}
\@ifundefined{tsGraphicsScaleY}{\gdef\tsGraphicsScaleY{.9}}{}
\def\checkGraphicsWidth{\ifdim\Gin@nat@width>\linewidth
	\tsGraphicsScaleX\linewidth\else\Gin@nat@width\fi}

\def\checkGraphicsHeight{\ifdim\Gin@nat@height>.9\textheight
	\tsGraphicsScaleY\textheight\else\Gin@nat@height\fi}

\def\fixFloatSize#1{}
\let\ts@includegraphics\includegraphics

\def\inlinegraphic[#1]#2{{\edef\@tempa{#1}\edef\baseline@shift{\ifx\@tempa\@empty0\else#1\fi}\edef\tempZ{\the\numexpr(\numexpr(\baseline@shift*\f@size/100))}\protect\raisebox{\tempZ pt}{\ts@includegraphics{#2}}}}

\AtBeginDocument{\def\includegraphics{\@ifnextchar[{\ts@includegraphics}{\ts@includegraphics[width=\checkGraphicsWidth,height=\checkGraphicsHeight,keepaspectratio]}}}

\DeclareMathAlphabet{\mathpzc}{OT1}{pzc}{m}{it}

\def\URL#1#2{\@ifundefined{href}{#2}{\href{#1}{#2}}}

\def\UrlOrds{\do\*\do\-\do\~\do\'\do\"\do\-}%
\g@addto@macro{\UrlBreaks}{\UrlOrds}

\edef\fntEncoding{\f@encoding}

\makeatother

\newif\ifmultipleabstract\multipleabstractfalse%
\usepackage{amsfonts,amssymb,amsmath}
\usepackage{mathtools,float}
\usepackage{graphicx}
\usepackage{caption}
\usepackage{amsthm}
\usepackage{enumitem}
\newtheorem{definition}{Definition}
\newtheorem{theorem}{Theorem}
\newtheorem{lemma}{Lemma}
\newtheorem{corollary}{Corollary}
\makeatletter
\newcommand*{\rom}[1]{\expandafter\@slowromancap\romannumeral #1@}
\makeatother
\newenvironment{myproof}[1][\proofname]{%
	\proof[\normalfont \bfseries #1]%
}{\endproof}

\begin{document}

\begin{frontmatter}

    \title{
  Improved Lower Bounds on the Domination Number of Hypercubes and Binary Codes with Covering Radius One    
}
    
\author[]{Ying-Sian Wu}
\ead{b11201018@ntu.edu.tw}
\ead{puggywu15@gmail.com}
\author[]{Jun-Yo Chen}
\ead{aggacc@gmail.com}

\begin{abstract}
 A dominating set on an $n $-dimensional hypercube is equivalent to a binary covering code of length $n $ and covering radius 1. It is still an open problem to determine the domination number $\gamma(Q_n) \text{ for } n\geq10 \text{ and } n\ne2^{k},2^{k}-1 $ ($k\in\mathbb{N} $). When $n$ is a multiple of 6, the best known lower bound is $\gamma(Q_n)\geq \frac{2^n}{n}$, given by Van Wee  (1988). In this article, we present a new method using congruence properties due to Laurent Habsieger (1997) and obtain an improved lower bound $\gamma(Q_n)\geq \frac{(n-2)2^n }{n^2-2n-2}$ when $n$ is a multiple of 6.
\end{abstract}
  \end{frontmatter}

\section{Introduction}
	Determining the domination number is an important optimization problem in graph theory, as well as an NP-complete problem in computational complexity theory \cite{David}. The domination problem on hypercubes is equivalent to the covering code problem. Generic introductions to domination problems and covering codes can be found in \cite{Haynes}\cite{CohenBook}.
	
	\textit{The} $n$\textit{-dimensional hypercube} $Q_n$ is defined recursively in terms of the cartesian product of graphs as follows,
	\begin{equation}
		Q_1=K_2,\;\;\; Q_n=K_2\square Q_{n-1}. \tag{1.1}
	\end{equation}
	Therefore, $Q_n$ can also be defined as
	\begin{equation}
		V(Q_n)=2^{\left\{1,2,\ldots,n\right\}}, \text{ }
		E(Q_n)=\left\{uv:u,v\in V(Q_n), u\subset v, \text{ and }	 |v\setminus u|=1\right\}.\tag{1.2}
	\end{equation}
	To avoid confusion, we use small brackets to express vertices. For instance, the vertex $\left\{2,3,5\right\}$ is written as $(2,3,5)$. Note that the vertex $\emptyset$ is written as $(0)$. Moreover, we define $(a_1,\ldots,a_i,0)\equiv(a_1,\ldots,a_i)$ to simplify some of our arguments.
	
	Given $S\subseteq V(Q_n)$, we define the function $g$ to express all different members in the union of the coordinates of the vertices in S. That is, 
	\begin{equation}
		g(S)\coloneqq\bigcup_{v\in S}v.\tag{1.3}
	\end{equation}
	Moreover, let $S[a]$ indicate the subset of S in which the number $a$ is contained in the coordinates of every vertex. That is,
	\begin{equation}
		S[a]\coloneqq \left\{ v : a\in v\text{ and } v\in S\right\}.\tag{1.4}
	\end{equation} 
	e.g. For $S=\Bigl\{ (1,2,3),(2,5),(3,5)\Bigr\}$, $g(S)=\left\{1,2,3,5\right\}$ and $S[5]=\Bigl\{(2,5),(3,5)\Bigr\}$.
	
	For some $v\in V(Q_n)$ and $S\subseteq V(Q_n)$, we define the neighborhood of $v$ and $S$ as follows.
	\begin{equation}
		N_i[v]=\left\{u\in V(Q_n): d(u,v)=i\right\}, N_i[S]=\bigcup_{v\in S}N_i[v].\notag
	\end{equation}	
	\begin{equation}
		N[v]=N_0[v]\cup N_1[v], N[S]=\bigcup_{v\in S}N[v].\tag{1.5}
	\end{equation}
	Given $S\subseteq V(Q_n)$, if $N[S]=V(Q_n)$, then we call $S$ a $dominating$ $set$ of $Q_n$.
	If there does not exist $ S'\subseteq V(Q_n)$ such that $|S'|<|S|$ and $N[S']=V(Q_n)$, then we call $S$ a $minimum$ $dominating$ $set$ of $Q_n$, and call $|S|$ the $domination$ $number$ of $Q_n$. In this article, we denote a given dominating set as $D$, and denote the domination number as $\gamma(Q_n)$.
	Tables 1 and 2 in the appendix give the latest results of the upper and lower bounds on $\gamma(Q_n)$. The summarization is due to Gerzson Kéri \cite{Keritable}\cite{keribook}.
	
	\begin{definition}
		Given a dominating set $D$ and $S\subseteq V(Q_n)$, we denote the $excess$ of $D$ on $S$ by $\delta_S(D)$, which is defined as $\sum_{v\in S}(|N[v]\cap D|-1)$. When $D$ is clear, we briefly write $\delta_S\coloneqq\delta_S(D)$. Also, if $S=\left\{v\right\}$, then $\delta_v\coloneqq\delta_{\left\{v\right\}}$.
	\end{definition}
	The term $excess$ has been used in many previous works, such as \cite{Laurent}\cite{gerard}.	We further define the symbols below. Likewise, the $D$ can be omitted if it's clear.
	\begin{equation}
		V\delta^x(D)\coloneqq\left\{ v: \delta_{v}(D)=x \text{ and } v\in V(Q_n) \right\},\notag
	\end{equation}
	\begin{equation}
		V\delta(D)\coloneqq \bigcup_{x\geq1}V\delta^x(D), \;\; C(D)\coloneqq \bigcup_{x\geq2}V\delta^x(D). \tag{1.6}
	\end{equation}
	\\
	Previous studies came up with various congruence properties of $\delta_{N_i[v]}(D)$, which help to estimate $\delta_{V(Q_n)}(D)$, and thus obtained the lower bounds on $\gamma(Q_n)$ due to the relation $\delta_{V(Q_n)}(D)=(n+1)|D|-|V(Q_n)|$ from \cite{gerard}. Theorem 1 is a segment of the properties given by Laurent Habsieger \cite{Laurent}, which we will apply.
	
	\begin{theorem}\upshape\textbf{(Habsieger)}\mbox{ }\itshape
		When $n$ is a multiple of $6$,
		\begin{equation}
			\delta_{N[v]}(D)   \equiv 1 \text{ }(mod\text{ } 2)  \text{, if } v\notin D.
			\tag{1.7}
		\end{equation}
		\begin{equation}
			\delta_{N[v]}(D)  \equiv0 \text{ }(mod\text{ } 2)  \text{, if } v\in D.
			\tag{1.8}
		\end{equation}
		\begin{equation}
			\delta_{N_1[v]}(D)+\delta_{N_2[v]}(D)  \equiv0 \text{ }(mod\text{ } 3).
			\tag{1.9}
		\end{equation}
	\end{theorem}
	
		We then put forward the pivotal concept throughout this article, $surfeit$. Although it seems closely related to $excess$, we shall demonstrate that such further analyzation is enough to improve the known bounds.
		
	\begin{definition}
		Given a dominating set $D$ and $S\subseteq V(Q_n)$, we denote the $surfeit$ of $D$ on $S$ by $\zeta_S(D)$, which is defined as $\sum_{v\in S\setminus D}(\delta_{N[v]}(D)-1)$. When D is clear, we briefly write $\zeta_S\coloneqq \zeta_S(D)$.
	\end{definition}
	
	 We further define the symbols below. Likewise, the $D$ can be omitted if it's clear.
	\begin{equation}
		V\zeta^x(D)\coloneqq\left\{v:\delta_{N[v]}(D)=x+1 \text{ and } v\in V(Q_n)\setminus D\right\},\text{ }
		V\zeta(D)\coloneqq \bigcup_{x\geq1}V\zeta^x(D) \tag{1.10}
	\end{equation}
	
	We look into the cases when $n$ is a multiple of 6. By calculating $\zeta_{V(Q_n)}(D)$ using two different methods, we show that it leads to a contradiction if $\gamma(Q_n)$ is too small.
	
\section{Generalities}
	
	All the arguments are considered in the cases when $n$ is a multiple of 6.
	\\
	
	Given an arbitrary value $\gamma^*$, we assume that there exists a dominating set $D$ satisfying $|D|=\gamma^*$. We calculate $\zeta_{V(Q_n)}$ using two different methods, and write the value we obtain as $\zeta_{m1}$ and $\zeta_{m2}$, respectively. A dominating set should lead to $\zeta_{m1}=\zeta_{m2}$. However, we will show that there must be $\zeta_{m2}>\zeta_{m1}$ when $\gamma^*$ is too small, implying that such dominating set $D$ cannot exist, so $\gamma(Q_n)>\gamma^*$.
	\\
	
	For all $ v\in V(Q_n)\setminus D,$ we have $\delta_{N[v]} \geq1$ by (1.7), so the first method to calculate $\zeta_{V(Q_n)}$ holds. Note that for all $ u\in V\delta^x,$ we have $|N[u]\setminus D|=n-x$, and $\sum_{x\in\mathbb{N}}x|V\delta^x|=\delta_{V(Q_n)}$. The value obtained this way is written as $\zeta_{m1}$:
	\begin{equation}
		\begin{aligned}
			\zeta_{V(Q_n)}&=\sum_{x\in\mathbb{N}}\sum_{u\in V\delta^x}x(n-x)-|V(Q_n)\setminus D|\\&=\sum_{x\in\mathbb{N}}x(n-x)|V\delta^x|-2^n+|D|\\&\noindent=(n-1)\delta_{V(Q_n)}-2^n+|D|-\sum_{x\in\mathbb{N}}x(x-1)|V\delta^x|\eqqcolon\zeta_{m1}.
		\end{aligned}
		\tag{2.1}
	\end{equation} 
	\\
	$\zeta_{m1}$ attains its maximum when $C=\emptyset$. We write this value as $\zeta_{\text{\text{max}}}$.
	\begin{equation}
		\zeta_{m1}\leq\zeta_{\text{\text{max}}}\coloneqq (n-1)\delta_{V(Q_n)}-2^n+|D|.\tag{2.2}
	\end{equation}
	
	Let us consider another method to calculate $\zeta_{V(Q_n)}$. The value obtained this way is written as $\zeta_{m2}$:
	\begin{equation}
		\zeta_{V(Q_n)}=\sum_{i\geq0}\Bigl(2i|V\zeta^{2i}|+(2i+1)|V\zeta^{2i+1}|\Bigr)=\sum_{i\geq1}2i|V\zeta^{2i}|\eqqcolon\zeta_{m2}.\tag{2.3}
	\end{equation}
	Note that by (1.7) we have $V\zeta^{2i+1}=\emptyset$.

	\begin{lemma}
		For all $ u\in V\delta^1$, if $d(u,C)\geq3$, then $|N[u]\cap V\zeta|\geq3$. 
	\end{lemma}
	\begin{myproof}
		Assume without loss of generality that $u=(0)$, $N[u]\cap D=\left\{(1),(a)\right\}$, where $a\in\left\{0,2,3,\ldots,n\right\}$, then $\left\{(0),(1,a)\right\}\subset V\delta^1$. For convenience we assume that $a\ne 0$, since the case $a=0$ can be dealt with similarly.
		Let $S= V\delta\cap(N_1[u]\cup N_2[u]) \setminus \left\{(1,a)\right\}$. 
		Applying (1.7) and (1.8) on the vertices in $N_1[u]$, we have the following:
		\begin{equation}
			\text{ For all }k\in\left\{1,2,\ldots,n\right\}\text{, there is } |S[k]|\equiv 0\text{ (mod 2)}.\tag{2.4}
		\end{equation}
		Moreover, by applying (1.9) on $u$, we have $|S|\equiv 2$ (mod 3).  
		In particular, $|S|\geq 5$, otherwise there exists $ k\in \left\{1,2,\ldots,n\right\}$ such that $|S[k]|=1$, contradicting (2.4).
		
		Suppose that $|N[u]\cap V\zeta|\leq2$, $N[u]\cap V\zeta\subseteq\left\{(b),(c)\right\}$ where $b,c\in\left\{0,2,3,\ldots,n\right\}\setminus\left\{(a)\right\}$, then $g(S)\subseteq\left\{1,a,b,c\right\}$. So let $T=\left\{(1),(a),(b),(c),(1,b),(1,c),(a,b),(a,c),(b,c)\right\}$, then $S\subseteq T$. If $b=0$, then $S=T\setminus\left\{(0)\right\}=\left\{(1),(a),(c),(1,c),(a,c)\right\}$, contradicting (2.4). Hence $b,c\ne 0$, $(0)\notin N[u]\cap V\zeta$, $|S\cap N_1[u]|=0$, $S=T\cap N_2[u]=\left\{(1,b),(1,c),(a,b),(a,c),(b,c)\right\}$, but this still contradicts (2.4).
		Therefore, $|N[u]\cap V\zeta|\geq3$.
	\end{myproof}
	In other words, for each $u\in V\delta^1$, there must be $d(u,C)\leq 2$ if $|N[u]\cap V\zeta|\leq2$. In Lemma 2 we will show that $|N[u]\cap V\zeta|\leq2$ gives us more rigorous conditions. To simplify our arguments, given $v\in C$, we divide $(N_1[v]\cup N_2[v])\cap V\delta^1$ into the following vertex sets: 

	$$T_1(v)\coloneqq(N_1[v]\setminus D)\cap V\delta^1;$$
	$$T_2(v)\coloneqq\left\{u\in N_2[v]\cap V\delta^1: |N[u]\cap N[v]\cap D|=0 \right\};$$
	\begin{equation}
		T_3(v)\coloneqq\left\{u\in N_2[v]\cap V\delta^1: |N[u]\cap N[v]\cap D|=2 \right\};\tag{2.5}
	\end{equation}
	$$T_4(v)\coloneqq\left\{u\in N_2[v]\cap V\delta^1: |N[u]\cap N[v]\cap D|=1 \right\};\text{ and}$$
	$$T_5(v)\coloneqq(N_1[v]\cap D)\cap V\delta^1.$$
	\\
	\begin{lemma}
		For all $ u\in V\delta^1$, if $|N[u]\cap V\zeta|\leq2$, then $|N[u]\cap V\zeta|=2$, and the following four claims hold.
		\begin{enumerate}
			\item There exists $ v\in C\setminus D$ such that $u\in T_1(v)\cup T_2(v)$, or there exists $ v\in C\cap D$ such that $u\in T_2(v).$
			\item For all $ v\in C\setminus D$ such that $u\in T_1(v)$, we have
			\begin{equation}
				|N_1[v]\cap V\delta|\leq3.\tag{2.6}
			\end{equation}
			\item For all $ v\in C$ such that $u\in T_2(v)\setminus D$, if we rename the coordinates so that $v=(0),$ $u=(a,b)$, then
			\begin{equation}
				(a),(b)\notin V\delta\text{ and }|(N_2[v]\cap V\delta)[a]|,|(N_2[v]\cap V\delta)[b]|\leq3.\tag{2.7}
			\end{equation}
			\item For all $ v\in C$ such that $u\in T_2(v)\cap D$, if we rename the coordinates so that $v=(0),$ $u=(a,b)$, then
			\begin{equation}
				|(N_2[v]\cap V\delta)[a]|,|(N_2[v]\cap V\delta)[b]|\leq2.\tag{2.8}
			\end{equation}
		\end{enumerate}
	\end{lemma}
	
	\newpage
	
	\begin{myproof}
		Given $v\in C$, we will prove the following statements.
		\begin{enumerate}[label=(\Alph*)]
			\item If $u\in T_1(v)$, $v\notin D$ and $|N[u]\cap V\zeta|\leq2$, then $|N[u]\cap V\zeta|=2$ and $(2.6)$ holds for $v$.
			\item If $u\in T_2(v)\setminus D$ and $|N[u]\cap V\zeta|\leq2$, then $|N[u]\cap V\zeta|=2$ and $(2.7)$ holds for $v$.
			\item If $u\in T_2(v)\cap D$ and $|N[u]\cap V\zeta|\leq2$, then $|N[u]\cap V\zeta|=2$ and $(2.8)$ holds for $v$.
			\item If $u\in T_1(v)$, $v\in D$ and $|N[u]\cap V\zeta|\leq2$, then $|N[u]\cap V\zeta|=2$ and there exists some $v'\in C\setminus D$ such that $u\in T_1(v')$.
			\item If $u\in T_3(v)\cup T_4(v)\cup T_5(v)$ and $u\notin T_1(v')\cup T_2(v')$ for all $v'\in C$, then $|N[u]\cap V\zeta|\geq3$.
		\end{enumerate}
	    (A), (B), and (C) prove Claims 2, 3, 4 in Lemma 2. By Lemma 1 there exists $v_0\in C$ with $u\in \bigcup_{1\leq i\leq 5}T_i(v_0)$, so Claim 1 follows by an application of (D) and (E) with $v=v_0$.\\[3pt] 
	    Below, (A) and (D) follow from Case 1, (B) and (C) follow from Case 2, while (E) follows from Cases 3 to 5.
	    We assume without loss of generality that $v=(0)$. 
		
		\textbf{\\Case 1-(1)} $u\in T_1(v)$ and $v\notin D$.
		\\
		
		$ \left\{u,v\right\}\subseteq N[u]\cap V\zeta$, so $|N[u]\cap V\zeta|\geq 2$.
		If equality holds, then $ N[u]\cap V\zeta=\left\{u,v\right\}$. Let $u=(a)\notin D$.
		Since $u\in V\delta^1$, we have $|g(N[u]\cap D)|=3$. So if $|N_1[v]\cap V\delta|\geq4$, then there exists $(k)\in N_1[v]\cap V\delta$ such that $k\notin g(N[u]\cap D)\text{, }(a,k)\in  N[u]\cap V\zeta$, which is a contradiction. Hence $|N_1[v]\cap V\delta|\leq 3$ and (A) is proved.
		
		\textbf{\\Case 1-(2)} $u\in T_1(v)$ and $v\in D$.
		\\
		
		If $\delta_v\geq 3$, then assume without loss of generality that $(1),(2),(3)\in D$ and that $u=(4)$. $u\notin C$, so $|\left\{(1,4),(2,4),(3,4)\right\}\cap D|\leq 1$. Note that $\left\{u,(1,4),(2,4),(3,4)\right\}\setminus D\subseteq N[u]\cap V\zeta$, so $| N[u]\cap V\zeta|\geq3$. Therefore, if $| N[u]\cap V\zeta|\leq2$, then $\delta_v=2$. 
		
		Let $N[v]\cap D=\left\{(0),(1),(2)\right\}$ and $u=(a)\notin D$.
		We have $(1),(2),(a)\in V\delta$, so for $w\in \left\{(1,a),(2,a)\right\}$, there is $w\notin  N[u]\cap V\zeta$ if and only if $ w\in D$. Thus, we assume that $\left\{(1,a),(2,a)\right\}\cap D=\left\{(1,a)\right\}$, $N[u]\cap V\zeta=\left\{u,(2,a)\right\}$.

		There does not exist $k \in \left\{3,4,\ldots,n\right\}\setminus\left\{a\right\}$ such that $(a,k)\in N[u]\cap V\zeta$. Therefore,
		\begin{equation}
			(N_1[u]\cup N_2[u])\cap V\delta\subseteq \left\{(0),(1),(2),(1,a),(2,a),(1,2,a)\right\}.\tag{*}
		\end{equation}
		Applying (1.7) and (1.9) on $u$, we get
		\begin{equation}
			\delta_{(0)}+\delta_{(1)}+\delta_{(2)}+\delta_{(1,a)}+\delta_{(2,a)}+\delta_{(1,2,a)}\equiv 0 \text{ (mod 3)},\tag{2.9}
		\end{equation}
		and 
		\begin{equation}
			\delta_{(1,a)}+\delta_{(2,a)}\equiv 0\text{ (mod 2). }\tag{2.10}
		\end{equation}
		\\
		 Moreover, $1\leq \delta_{(1,a)}\leq 2$ since by (*) we know that
		$N[(1,a)]\cap D \subseteq\left\{(1),(1,a),(1,2,a)\right\}$.
		\newpage 
		
		If $\delta_{(1,a)}=2$, then $(1,2,a)\in D$ and thus $\delta_{(2,a)}\geq1$. By (2.10) we know $\delta_{(2,a)}\geq2$. 
		
		If $\delta_{(1,a)}=1$ and $(1,2)\notin D$, then by (*) we know $\delta_{(1)}=2$, $\delta_{(2)}=1$, and $\delta_{(1,2,a)}=0$. Now that $\delta_{(0)}+\delta_{(1)}+\delta_{(2)}+\delta_{(1,a)}+\delta_{(1,2,a)}=6$, and by (2.10) we know $\delta_{(2,a)}\geq 1$, so (2.9) suggests that $\delta_{(2,a)}\geq 3$.
		
		If  $\delta_{(1,a)}=1$ and $(1,2)\in D$, then by (*) we know $\delta_{(1)}=3$, $\delta_{(2)}=2$, and $\delta_{(1,2,a)}=1$. Now that $\delta_{(0)}+\delta_{(1)}+\delta_{(2)}+\delta_{(1,a)}+\delta_{(1,2,a)}=9$, and by (2.10) we know $\delta_{(2,a)}\geq 1$, so (2.9) suggests that $\delta_{(2,a)}\geq 3$.
	
	    Therefore, there must be $(2,a)\in C\setminus D$ and $u\in T_1((2,a))$. This shows that there exists $v'\in C\setminus D$ such that $u\in T_1(v')$, and (D) is proved.
		
		\textbf{\\Case 2.} $u\in T_2(v)$.
		\\
		
		Let $u=(a,b)$ where $(a),(b)\notin D$, then $\left\{(a),(b)\right\}\subseteq N[u]\cap V\zeta$, so $| N[u]\cap V\zeta|\geq2$. Assume that equality holds. 
		Consider the case $u\notin D$.
		We have $(a),(b)\notin V\delta$, for otherwise $u\in  N[u]\cap V\zeta$, which is a contradiction. Also,
		$|g(N[u]\cap D)|=4$, so if $|(N_2[v]\cap V\delta)[a]|\geq4$, then there exists $ w\in (N_2[v]\cap V\delta)[a]$, $k\in w\setminus\left\{a\right\}$ such that
		\begin{equation}
			k\notin g(N[u]\cap D)\text{, }(a,b,k)\in  N[u]\cap V\zeta,\notag
		\end{equation}
		which is a contradiction. Likewise, we have $|(N_2[v]\cap V\delta)[b]|\leq3$, so
		(B) is proved. The same argument can be applied to the case $u\in D$ and prove (C).
		\\
		
		The following cases together prove (E).
		
		\textbf{\\Case 3.} $u\in T_3(v)$ and $u \notin T_1(v')\cup T_2(v')$ for all $v'\in C$.
		\\
		
		Let $(1),(2)\in D\text{ and } u=(1,2)$. We prove $| N[u]\cap V\zeta|\geq3$ by contradiction, assuming $| N[u]\cap V\zeta|\leq2$.  If $(a)\in D$ for some $a\in \left\{0,3,4,\ldots,n\right\}$, then $(1,2,a)\in N[u]\cap V\zeta$. This implies $\delta_v\leq3$, and since such an $a$ exists by $\delta_v\geq 2$, we have $| N[u]\cap V\zeta|\geq1$.
		
		If $| N[u]\cap V\zeta|=1$, we let $ N[u]\cap V\zeta=\left\{(1,2,a)\right\}$, where $a\in\left\{0,3,4,\ldots,n\right\}$, then $(N_1[u]\cup N_2[u])\cap V\delta\subseteq\left\{(0),(1,a),(2,a)\right\}$ and $\delta_{(1,a)},\delta_{(2,a)}\geq1$. We know $\left\{(1,a),(2,a)\right\}\cap C\ne\emptyset$ by applying (1.9) on $u$. Let $(1,a)\in C$, then there exists
		$ b\in \left\{0,3,4,\ldots,n\right\}\setminus\left\{a\right\}$ such that
		\begin{equation}
			(1,a,b)\in D\text{, }(1,b)\in V\delta\text{, }(1,2,b)\in  N[u]\cap V\zeta,\notag
		\end{equation}
		which is a contradiction. Therefore, $| N[u]\cap V\zeta|=2$.\\
		
		Now let $ N[u]\cap V\zeta=\left\{(1,2,3),(1,2,k)\right\}$, where $k\in\left\{0,4,5,\ldots,n\right\}$, then
		\begin{equation}
			(N_1[u]\cup N_2[u])\cap V\delta\subseteq\left\{(0),(1,3),(2,3),(1,k),(2,k),(1,2,3,k)\right\}.\tag{**}
		\end{equation}
		Thus, $N[v]\cap D\subseteq\left\{(1),(2),(3),(k)\right\}$. We denote the excesses as $\delta_{(0)}\eqqcolon o$, $\delta_{(1,3)}\eqqcolon p$, $\delta_{(2,3)}\eqqcolon q$, $\delta_{(1,k)}\eqqcolon r$, $\delta_{(2,k)}\eqqcolon s$, $\delta_{(1,2,3,k)}\eqqcolon t$, respectively. By (1.7), (1.8) and (1.9) we derive the following relations (note that $\delta_{(1,2)}=1$):
		\begin{align}
			o+p+r\equiv 1\text{ (mod 2) }\text{, }&\text{ since }\delta_{(1)}+\delta_{N_1[(1)]} \equiv 0\text{ (mod 2) };\tag{\rom{1}}\\
			o+q+s\equiv 1\text{ (mod 2) }\text{, }&\text{ since }\delta_{(2)}+\delta_{N_1[(2)]} \equiv 0\text{ (mod 2) };\tag{\rom{2}}\\
			p+q+t\equiv 0\text{ (mod 2) }\text{, }&\text{ since }\delta_{(1,2,3)}+\delta_{N_1[(1,2,3)]} \equiv 1\text{ (mod 2) };\tag{\rom{3}}\\
			r+s+t\equiv 0\text{ (mod 2) }\text{, }&\text{ since }\delta_{(1,2,k)}+\delta_{N_1[(1,2,k)]} \equiv 1\text{ (mod 2) };\tag{\rom{4}}\\
			o+p+q+r+s+t\equiv 0& \text{ (mod 3), }\text{ since }\delta_{N_1[u]}+\delta_{N_2[u]}\equiv 0\text{ (mod 3) }.\tag{\rom{5}}
		\end{align}
		
		If $o=3$, then $N[v]\cap D=\left\{(1),(2),(3),(k)\right\}$, and the following relations hold. Note that the restrictions leading to these results are due to (**).
		\begin{align}
			1\leq p\leq2\text{, }&\text{ since }N[(1,3)]\cap D\subseteq \left\{(1),(3),(1,3,k)\right\}\text{ and }(1),(3)\in D;\notag\\
			1\leq q\leq2\text{, }&\text{ since }N[(2,3)]\cap D\subseteq \left\{(2),(3),(2,3,k)\right\}\text{ and }(2),(3)\in D;\notag\\
			1\leq r\leq2\text{, }&\text{ since }N[(1,k)]\cap D\subseteq \left\{(1),(k),(1,3,k)\right\}\text{ and }(1),(k)\in D;\notag\\
			1\leq s\leq2\text{, }&\text{ since }N[(2,k)]\cap D\subseteq \left\{(2),(k),(2,3,k)\right\}\text{ and }(2),(k)\in D;\notag\\
			t\leq 1\text{, }&\text{ otherwise }u\in T_1((1,2,3,k))\cup T_2((1,2,3,k)).\notag
		\end{align}
		
		Now if $o+p+q+r+s+t=12$, then $p=q=r=s=2$ and $t=1$, contradicting (\rom{3}), so the only possibility is $o+p+q+r+s+t=9$. However, this implies that $(p+q+t)+(r+s+t)-t=6$. By (\rom{3}) and (\rom{4}) we know $t=0$, and together with (\rom{1}) and (\rom{2}) we know that $p,q,r,s$ have the same parity, which is impossible. Therefore,
		$|N[u]\cap V\zeta|\geq3$.
		\\
		
		If $o=2$, then we assume without loss of generality that $N[v]\cap D \subseteq \left\{(1),(2),(k)\right\}$. This time we obtain $p\leq 1$, $q\leq 1$, $1\leq r\leq 2$, $1\leq s\leq 2$, $t\leq1$. If $o+p+q+r+s+t=9$, then $p=q=1$, $r=s=2$, $t=1$, contradicting (\rom{3}), so the only possibility left is $o+p+q+r+s+t=6$. This implies that $(p+q+t)+(r+s+t)-t=4$. By (\rom{3}) and (\rom{4}) we know $t=0$, and $p,q$ as well as $r,s$ have the same parity. Thus $p=q$ and $r=s$, implying $p+r=2$, which contradicts (\rom{1}). Therefore, $| N[u]\cap V\zeta|\geq3$.
		
		\textbf{\\Case 4.} $u\in T_5(v)$ and $u \notin T_1(v')\cup T_2(v')$ for all $v'\in C$. 
		\\
		
		If $v\in D$, then we can prove $| N[u]\cap V\zeta|\geq3$ using the same method as in Case 3. We sketch our arguments in a simplified version, for they are highly similar to those in Case 3:\\
		
		Let $(1),(2)\in D$ and $u=(1)$. We prove our claim by contradiction, assuming that $| N[u]\cap V\zeta|\leq2$. Like in case 3 we see $\delta_v\leq 3$ and $| N[u]\cap V\zeta|=2$. Assume without loss of generality that $N[u]\cap V\zeta=\left\{(1,2),(1,3)\right\}$, then
		\begin{equation}
			(N_1[u]\cup N_2[u])\cap V\delta\subseteq\left\{(0),(2),(3),(1,2),(1,3),(1,2,3)\right\}.\tag{***}
		\end{equation}
		We denote the excesses as $\delta_{(0)}\eqqcolon o$, $\delta_{(2)}\eqqcolon p$, $\delta_{(3)}\eqqcolon q$, $\delta_{(1,2)}\eqqcolon r$, $\delta_{(1,3)}\eqqcolon s$, $\delta_{(1,2,3)}\eqqcolon t$, respectively. By (1.7), (1.8), and (1.9) we derive the following relations (note that $\delta_{(1)}=1$):
		\begin{align}
			o+r+s\equiv 1\text{ (mod 2) }\text{, }&\text{ since }\delta_{(1)}+\delta_{N_1[(1)]} \equiv 0\text{ (mod 2) };\tag{\rom{6}}\\
			o+p+q\equiv 1\text{ (mod 2) }\text{, }&\text{ since }\delta_{(0)}+\delta_{N_1[(0)]} \equiv 0\text{ (mod 2) };\tag{\rom{7}}\\
			p+r+t\equiv 0\text{ (mod 2) }\text{, }&\text{ since }\delta_{(1,2)}+\delta_{N_1[(1,2)]} \equiv 1\text{ (mod 2) };\tag{\rom{8}}\\
			q+s+t\equiv 0\text{ (mod 2) }\text{, }&\text{ since }\delta_{(1,3)}+\delta_{N_1[(1,3)]} \equiv 1\text{ (mod 2) };\tag{\rom{9}}\\
			o+p+q+r+s+t\equiv 0& \text{ (mod 3), }\text{ since }\delta_{N_1[u]}+\delta_{N_2[u]}\equiv 0\text{ (mod 3) }.\tag{\rom{10}}
		\end{align}
		
		If $o=3$, then $N[v]\cap D=\left\{(0),(1),(2),(3)\right\}$, and the following relations hold. Note that the restrictions leading to these results are due to (***).
		\begin{align}
			1\leq p\leq2\text{, }&\text{ since }N[(2)]\cap D\subseteq \left\{(0),(2),(2,3)\right\}\text{ and }(0),(2)\in D;\notag\\
			1\leq q\leq2\text{, }&\text{ since }N[(3)]\cap D\subseteq \left\{(0),(3),(2,3)\right\}\text{ and }(0),(3)\in D;\notag\\
			1\leq r\leq2\text{, }&\text{ since }N[(1,2)]\cap D\subseteq \left\{(1),(2),(1,2,3)\right\}\text{ and }(1),(2)\in D;\notag\\
			1\leq s\leq2\text{, }&\text{ since }N[(1,3)]\cap D\subseteq \left\{(1),(3),(1,2,3)\right\}\text{ and }(1),(3)\in D;\notag\\
			t\leq 1\text{, }&\text{ otherwise }u\in T_2((1,2,3)).\notag
		\end{align}
		We have $o+p+q+r+s+t=9$ or $12$. The latter contradicts (\rom{8}), while the former suggests that $(p+r+t)+(q+s+t)-t=6$, and using (\rom{6}) to (\rom{9}) we know $t=0$ and $p,q,r,s$ have the same parity, which is impossible.
		
		If $o=2$, then $N[v]\cap D \subseteq \left\{(0),(1),(2)\right\}$. This time we obtain $1\leq p\leq 2$, $q\leq 1$, $1\leq r\leq 2$, $s\leq 1$, $t\leq1$, so we have $o+p+q+r+s+t=6$ or $9$, but again we can easily lead to contradictions using (\rom{6}) to (\rom{10}). Therefore, $| N[u]\cap V\zeta|\geq3$.
		\\		
		
		If $v\notin D$, then let $u=(1)\in D$. 
		By $\delta_v\geq2$ we may assume $N[v]\cap D\supseteq\left\{(1),(2),(3)\right\}$. For $ w\in \left\{(1,2),(1,3)\right\}$, we have $w\notin  N[u]\cap V\zeta$ if and only if $ w\in D$. Moreover, $(0)\in  N[u]\cap V\zeta$. So if $|\left\{(1,2),(1,3)\right\}\cap D|=0$, then $| N[u]\cap V\zeta|\geq3$. If not, let $(1,2)\in D$, then $(1,2)\in C\cap D$ and $u\in T_5((1,2))$, implying $| N[u]\cap V\zeta|\geq3$.
		
		\textbf{\\Case 5.} $u\in T_4(v)$ and $u \notin T_1(v')\cup T_2(v')$ for all $v'\in C$.
		\\
		
		If $v\in D$, then we let $(1),(2)\in N[v]\cap D$, $u=(1,a)$, where $(a)\notin D$. 
		\\
		We have $\left\{(0),(1),(1,2),(1,a)\right\}\subset V\delta$. Therefore, if $\left\{(1,a),(1,2,a)\right\}\cap D=\emptyset$, then $\left\{(a),(1,a),(1,2,a)\right\}\subseteq N[u]\cap V\zeta$; if $(1,2,a)\in D$, then $(1,2)\in C$ and $u\in T_3((1,2))$; if $(1,a)\in D$, then $(1)\in C$ and $u\in T_5((1))$.
		
		On the other hand, if $v\notin D$, then let $(1),(2),(3)\in N[v]\cap D$, $u=(1,a)$, where $(a)\notin D$. 
		We have $\left\{(1,2),(1,3)\right\}\subset V\delta$ and $(a)\in  N[u]\cap V\zeta$. Thus, if $\left\{(1,2,a),(1,3,a)\right\}\cap D=\emptyset$, then $\left\{(a),(1,2,a),(1,3,a)\right\}\subseteq N[u]\cap V\zeta$; if $\left\{(1,2,a),(1,3,a)\right\}\cap D\ne\emptyset$, then there exists $ w\in \left\{(1,2),(1,3)\right\}\cap C\text{ such that }u\in T_3(w)$.
		
		By Case 3 and Case 4, every possible condition above implies that $|N[u]\cap V\zeta|\geq3$.
	\end{myproof}\mbox{}
	
	Given $v\in C$, we define $S_i(v)\coloneqq\left\{u\in T_i(v):|N[u]\cap V\zeta|=2\right\}$. Lemma 3 gives an upper bound for $|S_1(v)\cup S_2(v)|$ when $v\notin D$ and an upper bound for $|S_2(v)|$ when $v\in D$.
	
	\newpage
	
	\begin{lemma}
		For $ v\in C$, we have $|S_1(v)\cup S_2(v)|\leq \frac{3}{2}(n-\delta_v)$ if $v\notin D$, and $|S_2(v)|\leq \frac{3}{2}(n-\delta_v)$ if $v\in D$.
	\end{lemma}
	
	\begin{myproof}
		Assume without loss of generality that $v=(0)$. 
		\\
		
		If $v\notin D$, then let $\left\{(1),(2),\ldots,(\delta_v+1)\right\}\subset D$. 
		
		Define $A\coloneqq g(S_2(v)\cap D)$ and $B\coloneqq g(S_2(v))\setminus A$. By (2.7) and (2.8) we know that
		\begin{equation}
			\forall k\in A\text{, we have }|S_2(v)[k]|\leq 2\text{; } \forall k\in B\text{, we have }|S_2(v)[k]|\leq 3.\notag
		\end{equation}
		So $|S_2(v)|\leq \frac{1}{2}(2|A|+3|B|)$.\\
		
		By (2.7) we also know that for all $ k \in g(S_1(v))$, we have $k\notin g(S_2(v)\setminus D)$, $k\notin B$. Hence, $|B|\leq n-\delta_v-1-|S_1(v)|$. Since $|A|+|B|\leq n-\delta_v-1$, we have
		\begin{equation}
			\begin{aligned}
				|S_1(v)\cup S_2(v)|&\leq|S_1(v)|+\frac{1}{2}\Bigl(2|A|+3|B|\Bigr)\\
				&\leq|S_1(v)|+\frac{1}{2}\Bigl(2|S_1(v)|+3(n-\delta_v-1-|S_1(v)|)\Bigr)\\&=|S_1(v)|+\frac{1}{2}\Bigl(3n-3\delta_v-3-|S_1(v)|\Bigr).	
			\end{aligned}
			\notag
		\end{equation}
		By (2.6) we know that $|S_1(v)|\leq 3$, so $|S_1(v)\cup S_2(v)|\leq\frac{3}{2}(n-\delta_v)$.
		\\
		
		On the other hand, if $v\in D$, then let $\left\{(0),(1),(2),\ldots,(\delta_v)\right\}\subset D$.
		
		By (2.7), for all $ k\in \left\{\delta_v+1,\delta_v+2,\ldots,n\right\}$, we have $|S_2(v)[k]|\leq 3$. Therefore, 
		$$|S_2(v)|\leq\frac{3}{2}|\left\{\delta_v+1,\delta_v+2,\ldots,n\right\}|=\frac{3}{2}(n-\delta_v),$$
		and Lemma 3 is proved. 
	\end{myproof}\mbox{}
	
	By definition, $\sum_{i\geq1}(2i+1)|V\zeta^{2i}|=\sum_{x\in\mathbb{N}}x\sum_{u\in V\delta^x}|N[u]\cap V\zeta|$, and we can now estimate its lower bound using the results in Lemma 3.
	
	\begin{lemma}
		For $n\geq 12$, the following inequality holds:
		\begin{equation}
			\begin{aligned}
				&\sum_{i\geq1}(2i+1)|V\zeta^{2i}|\\
				&\geq3\delta_{V(Q_n)}-|V\delta^2|-4.5|V\delta^{n-3}|-(n+1)|V\delta^{n-2}|-(2n-0.5)|V\delta^{n-1}|-3n|V\delta^n|.
			\end{aligned}
			\notag
		\end{equation}
	\end{lemma}

	\newpage
	
	\begin{myproof}
		\mbox{}
		By Claim 1 in Lemma 2, we have $$\left\{u : |N[u]\cap V\zeta|-3<0, u\in V\delta^1\right\}= \bigcup_{v\in C\setminus D}(S_1(v)\cup S_2(v))\cup \bigcup_{v\in C\cap D}(S_2(v)).$$ 
		Therefore,
		\begin{equation}
			\begin{aligned}
				&\sum_{i\geq1}(2i+1)|V\zeta^{2i}|\\
				&=\sum_{x\in\mathbb{N}}x\sum_{u\in V\delta^x}|N[u]\cap V\zeta|\\
				&=3\delta_{V(Q_n)}+\sum_{x\in\mathbb{N}}x\sum_{u\in V\delta^x}\Bigl(|N[u]\cap V\zeta|-3\Bigr)\\
				&=3\delta_{V(Q_n)}+\sum_{v\in C}\delta_v\Bigl(|N[v]\cap V\zeta|-3\Bigr)+\sum_{u\in V\delta^1}\Bigl(|N[u]\cap V\zeta|-3\Bigr)\\
				&\geq3\delta_{V(Q_n)}+\sum_{v\in C}\delta_v(n-\delta_v-3)+\sum_{v\in C\setminus D}\sum_{u\in S_1(v)\cup S_2(v)}\Bigl(|N[u]\cap V\zeta|-3\Bigr)\\
				&\text{ }\text{ }\text{ }+\sum_{v\in C\cap D}\sum_{u\in S_2(v)}\Bigl(|N[u]\cap V\zeta|-3\Bigr)\\
				&=3\delta_{V(Q_n)}+\sum_{v\in C\setminus D}\Bigl(\delta_v(n-\delta_v-3)-|S_1(v)\cup S_2(v)|\Bigr)+\sum_{v\in C\cap D}\Bigl(\delta_v(n-\delta_v-3)-| S_2(v)|\Bigr)\\
				&\geq3\delta_{V(Q_n)}+\sum_{v\in C}\Bigl(\delta_v(n-\delta_v-3)-\frac{3}{2}(n-\delta_v)\Bigr).
			\end{aligned}
			\notag
		\end{equation}	
		Note that the last inequality is due to Lemma 3. A short calculation shows that for $n\geq 12$,
		\begin{equation}
			\delta_v(n-\delta_v-3)-\frac{3}{2}(n-\delta_v)\geq
			\begin{cases}
				0,\text{  if } 3\leq\delta_v\leq n-4\text{, or } \delta_v=2\text{ and }n\geq 18;\\
				-1,\text{  if } \delta_v=2\text{ and } n=12;\\
				-4.5,\text{  if } \delta_v=n-3;\\
				-n-1,\text{  if } \delta_v=n-2;\\
				-2n+0.5,\text{  if }\delta_v=n-1;\\
				-3n,\text{  if }\delta_v=n,
			\end{cases}
		\notag
		\end{equation}
	and Lemma 4 follows.
	\end{myproof}\mbox{}

    Finally, we can estimate $\zeta_{m2}-\zeta_{m1}$.
    \newpage
	
	\begin{lemma}
		 When $n\geq 12$, $\zeta_{m2}-\zeta_{m1}\geq2\delta_{V(Q_n)}-\zeta_{\text{max}}.$
	\end{lemma}
	\begin{myproof}
		By (2.3) and Lemma 4, we have
		\begin{equation}
		\begin{aligned}
			&\zeta_{m2}-2\delta_{V(Q_n)}\\
			&=\frac{1}{3}\left(\sum_{i\geq1}(6i|V\zeta^{2i}|)-6\delta_{V(Q_n)}\right)\\
			&=\frac{1}{3}\left(\sum_{i\geq1}(2i-2)|V\zeta^{2i}|+2\sum_{i\geq1}(2i+1)|V\zeta^{2i}|-6\delta_{V(Q_n)}\right)\\
			&\geq\frac{1}{3}\Bigl(\sum_{i\geq1}(2i-2)|V\zeta^{2i}|-2|V\delta^2|-9|V\delta^{n-3}|-(2n+2)|V\delta^{n-2}|\\
			&\text{ }\text{ }\text{ }\text{ }-(4n-1)|V\delta^{n-1}|-6n|V\delta^n|\Bigr).
		\end{aligned}
		\notag
		\end{equation}
		By (2.1) and (2.2) we have
		\begin{equation}
			\begin{aligned}
				&\zeta_{m1}-\zeta_{\text{max}}\\
				&=-\sum_{x\in\mathbb{N}}x(x-1)|V\delta^x|\\
				&=-2|V\delta^2|-6|V\delta^3|-\ldots-(n-3)(n-4)|V\delta^{n-3}|-(n-2)(n-3)|V\delta^{n-2}|\\
				&\text{ }\text{ }\text{ }\text{ }-(n-1)(n-2)|V\delta^{n-1}|-n(n-1)|V\delta^n|.
			\end{aligned}
			\notag
		\end{equation}
		Therefore,  $\zeta_{m2}-2\delta_{V(Q_n)}\geq\zeta_{m1}-\zeta_{\text{max}}$ when $n\geq 12$ and Lemma 5 follows.
	\end{myproof}\mbox{}
	
	\begin{theorem}\mbox{}
		If $n\equiv 0$ (mod 6), then $\gamma(Q_n)\geq\dfrac{(n-2)2^n}{n^2-2n-2}$.
	\end{theorem}

	\begin{myproof}
		Theorem 2 holds true for $n=6$ by $\gamma(Q_6)=12$, so assume $n\geq12$.
		Consider a minimum dominating set of $Q_n$. We have
		$$\delta_{V(Q_n)}=(n+1)\gamma(Q_n)-2^n,
		\text{ }\zeta_{\text{max}}=(n-1)\delta_{V(Q_n)}-2^n+\gamma(Q_n)=n^2\gamma(Q_n)-n2^n.$$
		By Lemma 5, $\zeta_{m2}-\zeta_{m1}\geq2\delta_{V(Q_n)}-\zeta_{\text{max}}$, so there must be $\zeta_{\text{max}}-2\delta_{V(Q_n)}\geq0$,
		$$(n^2-2n-2)\gamma(Q_n)-(n-2)2^n\geq0,\text{ }\gamma(Q_n)\geq\dfrac{(n-2)2^n}{n^2-2n-2}.$$
	\end{myproof}
	\begin{corollary}
		$\gamma(Q_{12})\geq348$, $\gamma(Q_{18})\geq14666$, $\gamma(Q_{24})\geq 701709$, $\gamma(Q_{30})\geq35876816$.
	\end{corollary}

\section{Conclusion}
	Tables 1 and 2 in the appendix are due to Gerzson Kéri \cite{Keritable}\cite{keribook}. When $n$ is a multiple of 6, previously the best known result was $\gamma(Q_n)\geq\frac{2^n}{n}$, given by van Wee \cite{gerard}. Our lower bound is higher, and several improvements are listed in Corollary 1.
	
\section*{Acknowledgements}
    We thank Gerard Jennhwa Chang and Chun-Ju Lai for encouragements and comments. We are grateful to Gerzson Kéri, Patric Östergård and Jan-Christoph Schlage-Puchta for their useful correspondence. We would like to thank Laurent Habsieger, whose fine paper inspired us a lot. Lastly, we would like to thank the reviewers and the editors, who helped us to improve the quality of this manuscript.

\section*{Appendix}
\begin{figure}[H]
	\centering
	\includegraphics[width=1\linewidth]{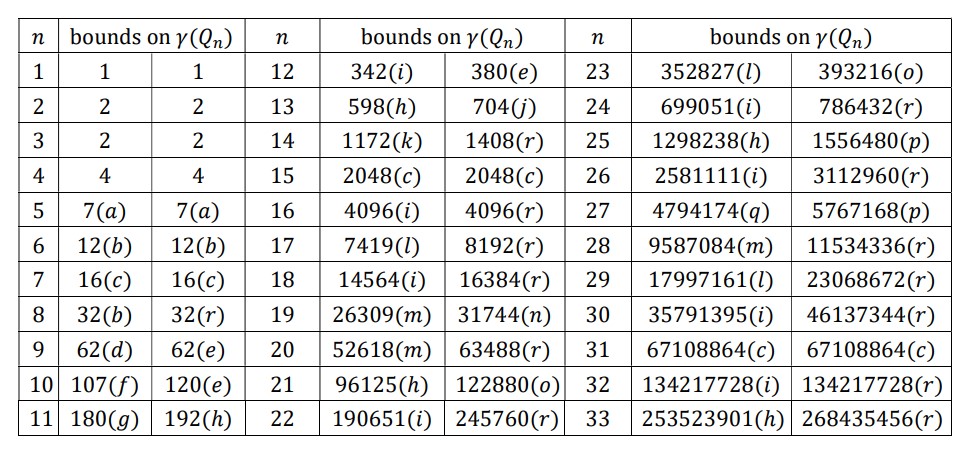}
	\captionof{table}[foo]{the latest results of the bounds on $\gamma(Q_n)$}
	\label{fig:uantwerpen-01}
\end{figure}

\begin{table}[H]
	\def\arraystretch{1.5}
	\centering
	\begin{tabular}{|c|c|c|c|}
		\hline
		$a$&Taussky, Todd, 1948 \cite{Todd}&$j$&Östergård, Weakly, 1999 \cite{Weakley2}\\ \hline
		$b$&Stanton, Kalbfleisch, 1968,1969 \cite{Stanton1}\cite{Stanton2}&$k$&Habsieger, 1997 \cite{Laurent}\\ \hline
		$c$&perfect code&$l$&Haas, 2007-2008 \cite{Haas}\\ \hline
		$d$&Östergård, Blass, 2001 \cite{Blass}&$m$&Habsieger, Plagne, 2000 \cite{Plagne}\\ \hline
		$e$&Wille, 1990, 1996 \cite{Wille1}\cite{Wille2}&$n$&Li, Chen, 1994 \cite{Li}\\ \hline
		$f$&Bertolo, Östergård, Weakley, 2004 \cite{Weakley}&$o$&Kéri, 2006 \cite{Keri2006}\\ \hline
		$g$&Blass, Litsyn, 1998 \cite{Litsyn}&$p$&Östergård, Kaikkonen, 1998 \cite{Kai}\\ \hline
		$h$&Cohen, Lobstein, Sloane, 1986 \cite{Sloane}&$q$&Plagne, 2009 \cite{PlagneR}\\ \hline
		$i$&van Wee,1988 \cite{gerard}&$r$&$\gamma(Q_{n+1})\leq 2\gamma(Q_n)$\\ \hline
	\end{tabular}
	\caption{References for Table 1}
\end{table}

\end{document}